\documentclass[12pt,reqno]{amsart}
\usepackage[margin=1in]{geometry}
\usepackage{amssymb,amsfonts,amsmath,mathrsfs,etoolbox,mathtools,bm}
\usepackage[breaklinks=true,colorlinks=true,linkcolor=teal,citecolor=teal,filecolor=teal,urlcolor=black!50!blue]{hyperref}
\hypersetup{linktocpage}
\usepackage[hiresbb]{graphicx,xcolor}
\allowdisplaybreaks[4]

\newtheorem{theorem}{Theorem}[section]
\newtheorem{lemma}[theorem]{Lemma}

\theoremstyle{definition}
\newtheorem{remark}[theorem]{Remark}

\numberwithin{equation}{section}

\patchcmd{\section}{\scshape}{\bfseries}{}{}
\makeatletter
\renewcommand{\@secnumfont}{\bfseries}
\makeatother
\newcommand{\norm}[1]{\left\lVert#1\right\rVert}
\makeatletter\newcommand{\tpmod}[1]{{\@displayfalse \pmod{#1}}}

\begin{document}

\title{On a Finite Field Analogue of van der Corput Differencing}

\author{Ikuya Kaneko}
\address{The Division of Physics, Mathematics and Astronomy, California Institute of Technology, 1200 E. California Blvd., Pasadena, CA 91125, USA}
\email{ikuyak@icloud.com}
\urladdr{\href{https://sites.google.com/view/ikuyakaneko/}{https://sites.google.com/view/ikuyakaneko/}}

\thanks{The author acknowledges the support of the Masason Foundation.}



\date{\today}

\dedicatory{}

\begin{abstract}
We establish nontrivial bounds for general bilinear forms with a given periodic function, which are thought of as an analogue of van der Corput differencing for exponential sums. The proof employs Poisson summation, Cauchy--Schwarz, and the iteration thereof.
\end{abstract}

\maketitle

\section{Introduction}\label{introduction}
The notion of exponential sums is one of the most important subjects in analytic number theory. The problem of determining the optimal size of such sums is key to understanding many salient problems such as equidistribution of various elements within domains, counting problems, the distribution of primes in subsets, and the number of solutions to Diophantine equations, to name a few. In modern number theory, there are various techniques to obtain their strong upper bounds; see in particular the pioneering papers of Weyl~\cite{Weyl1916} and van der Corput~\cite{vanderCorput1921}.~They consider a general exponential sum of the form
\begin{equation*}
\sum_{N \leq n < 2N} e(f(n)),
\end{equation*}
where $f(n)$ is a real-valued function and $e(x) = e^{2\pi ix}$. In order to prove a nontrivial estimate $N^{1-\delta}$ for some $\delta > 0$, van der Corput begins with replacing the phase~$f(x)$ with its derivative $f'(x)$ via the Cauchy--Schwarz inequality, and then applies the Poisson summation formula so that the dual length is roughly of size $f'(N)$. The iteration of van der Corput differencing leads to a standard $k$-th derivative test, which now evolves into the theory of exponent pairs. For further details and historical background, see the book of Graham--Kolesnik~\cite{GrahamKolesnik1991}.

This work addresses a finite field analogue of van der Corput differencing for exponential sums. To enable the ensuing discussion, we introduce some notation. For a sequence $\alpha = (\alpha_{n})$ supported on a dyadic interval $n \in [N, 2N)$, we write
\begin{equation*}
\norm{\alpha}_{2}^{2} \coloneqq \frac{1}{N} \sum_{N \leq n < 2N} |\alpha_{n}|^{2},
\end{equation*}
and for a $c$-periodic function $K(m, n; c)$ in $m$ and $n$, we write
\begin{equation}\label{eq:K-norm}
\norm{K}_{2, \infty}^{2} \coloneqq \sup_{(m, c) = 1} \frac{1}{c} \sum_{x \tpmod{c}} |K(m, x; c)|^{2}.
\end{equation}
Note that $\norm{\alpha}_{2}^{2}$ is normalised by a factor of $\frac{1}{N}$, which is less standard. Given such a function $K$, we denote
\begin{equation*}
(\Lambda_{\ell} K)(m, n; c) \coloneqq \frac{1}{\sqrt{c}} \sum_{x \tpmod{c}} K(m, x; c) \overline{K(n, x; c)} e \left(-\frac{\ell x}{c} \right).
\end{equation*}
Generally, one expects that if $K$ is a trace function, then so is $\Lambda_{\ell} K$ and thus $\Lambda_{\ell} K \ll 1$. We then define a bilinear form by
\begin{equation}\label{eq:bilinear-form}
\mathcal{S}_{K}(M, N; c) \coloneqq \sum_{M \leq m < 2M} \sum_{N \leq n < 2N} \alpha_{m} \beta_{n} K(m, n; c).
\end{equation}
The main result of this paper is the following.
\begin{theorem}\label{thm:main}
Let $K(m, n; c)$ be a $c$-periodic function in $m$ and $n$. Let $\alpha = (\alpha_{m})$ and~$\beta = (\beta_{n})$ traverse sequences supported on $m \in [M, 2M)$ and $n \in [N, 2N)$, respectively. Suppose that $M \asymp N$. Then we have that
\begin{multline}\label{SK}
\mathcal{S}_{K}(M, N; c) \ll \sqrt{cMN} \norm{\alpha}_{2} \norm{\beta}_{2} \sup_{|\ell_{1}|, \ldots, |\ell_{k}| \ll \frac{c}{N}} 
\bigg(\norm{K}_{2, \infty}+\norm{\Lambda_{\ell_{1}} K}_{2, \infty}^{\frac{1}{2}}+\cdots\\
 + \norm{\Lambda_{\ell_{k-1}} \cdots \Lambda_{\ell_{1}} K}_{2, \infty}^{\frac{1}{2^{k-1}}}+M^{\frac{1}{2^{k+1}}} \norm{\alpha}_{2}^{-\frac{1}{2^{k}}} 
\norm{\Lambda_{\ell_{k}} \cdots \Lambda_{\ell_{1}} K}_{2, \infty}^{\frac{1}{2^{k}}} \bigg).
\end{multline}
\end{theorem}

The proof of Theorem~\ref{thm:main} features Poisson summation, Cauchy--Schwarz, and the iteration \textit{ad infinitum} thereof. Despite the softness of our input, it has a square-root improvement over the trivial bound of size $MN$ (if we assume that all the norms are bounded). Furthermore, Theorem~\ref{thm:main} keeps the coefficients $(\alpha_{m})$ and $\beta_{n}$ as general as possible, which would facilitate potential arithmetic implications for the subconvexity problem for families of $L$-functions for instance. In particular, a version of the bilinear form~\eqref{eq:bilinear-form} often arises when one utilises the delta method to separate oscillations, in the terminology of Munshi~\cite{Munshi2018}. In this respect, an interesting case is when $K$ is a hyper-Kloosterman sum, in which case~\eqref{eq:bilinear-form} is studied~by Kowalski--Michel--Sawin~\cite{KowalskiMichelSawin2017} and Kerr et al.~\cite{KerrShparlinskiWuXi2023}. Note that these iterative ideas are present in~\cite{AggarwalLeungMunshi2022,HeathBrown1978-2,KanekoLeung2023-2}, but not written down in a general format.

\begin{remark}
The assumption $M \asymp N$ in Theorem~\ref{thm:main} is imposed to decrease the complexity of the argument, and in fact our strategy should be capable of dealing with a wider regime.
\end{remark}

\subsection*{Acknowledgements}
The author thanks Maksym Radziwi{\l}{\l} for suggesting underlying ideas in this paper, and Wing Hong Leung for productive discussions.

\section{Preliminaries}
For $n \in \mathbb{N}$ and an integrable function $w \colon \mathbb{R}^{n} \to \mathbb{C}$, we denote its Fourier transform by
\begin{equation*}
\widehat{w}(\xi) \coloneqq \int_{\mathbb{R}^{n}} w(t) e(-\langle t, \xi \rangle) dt,
\end{equation*}
where $\langle \cdot, \cdot \rangle$ denotes the standard inner product on $\mathbb{R}^{n}$. Furthermore, if $c \in \mathbb{N}$ and $K \colon \mathbb{Z} \to \mathbb{C}$ is a periodic function of period $c$, then its discrete Fourier transform $\widehat{K}$ is again the periodic function of period $c$ given by
\begin{equation*}
\widehat{K}(n) \coloneqq \sum_{a \tpmod{c}} K(a) e \left(-\frac{an}{c} \right).
\end{equation*}
We invoke a version of the Poisson summation formula with a $c$-periodic function involved.
\begin{lemma}[{Fouvry--Kowalski--Michel~\cite[Lemma~2.1]{FouvryKowalskiMichel2015}}]\label{Fouvry-Kowalski-Michel}
For any $c \in \mathbb{N}$, any $c$-periodic function $K$, and any even smooth function $V$ compactly supported on $\mathbb{R}$, we have that
\begin{equation*}
\sum_{n = 1}^{\infty} K(n) V(n) = \frac{1}{c} \sum_{n \in \mathbb{Z}} \widehat{K}(n) \widehat{V} \left(\frac{n}{c} \right).
\end{equation*}
\end{lemma}

\section{Proof of Theorem~\ref{thm:main}}
To eschew excessive complications, we can assume without loss of generality that $\alpha_{m}$ is such that $\alpha_{m} \ne 0$ implies $m$ is prime and $(m, c) = 1$. We further assume that $|K(m, n; c)| \leq 1$~for $(m, n) = 1$, and $|K(m, n; c)| \leq \sqrt{c}$ otherwise. Throughout this section, the notation $m \sim M$ serves as a shorthand for $M \leq m < 2M$.

First of all, the Cauchy--Schwarz inequality implies that there exists a compactly supported smooth function $V \in C_{c}^{\infty}([1, 2])$ such that
\begin{equation*}
|\mathcal{S}_{K}(M, N; c)| \leq N \norm{\beta}_{2} \left(\frac{1}{N} \sum_{n} 
\left|\sum_{m \sim M} \alpha_{m} K(m, n; c) \right|^{2} V \left(\frac{n}{N} \right) \right)^{\frac{1}{2}}.
\end{equation*}
Squaring out and applying Poisson summation (Lemma~\ref{Fouvry-Kowalski-Michel}) to the sum over $n$, we obtain
\begin{equation*}
c^{\frac{1}{4}} \sqrt{N} \norm{\beta}_{2} \left(\frac{N}{c} \sum_{\ell} \sum_{m_{1}, m_{2} \sim M} 
\alpha_{m_{1}} \alpha_{m_{2}} (\Lambda_{\ell} K)(m_{1}, m_{2}; c) \widehat{V} \left(\frac{\ell}{c/N} \right) \right)^{\frac{1}{2}}.
\end{equation*}
Repeated integration by parts ensures an arbitrary saving unless $|\ell| \ll \frac{c}{N}$, thus we arrive at
\begin{equation*}
c^{\frac{1}{4}} \sqrt{N} \norm{\beta}_{2} \left(\frac{N}{c} \sum_{|\ell| \ll \frac{c}{N}} \sum_{m_{1}, m_{2} \sim M} 
\alpha_{m_{1}} \alpha_{m_{2}} (\Lambda_{\ell} K)(m_{1}, m_{2}; c) \right)^{\frac{1}{2}}.
\end{equation*}
We decompose this expression according as $m_{1} = m_{2}$ or $m_{1} \ne m_{2}$. Upon bounding the sum over $\ell$ trivially, the contribution of the diagonal terms $m_{1} = m_{2}$ boils down to
\begin{equation*}
\sqrt{cMN} \norm{\beta}_{2} \left(\frac{1}{M} \sum_{m \sim M} |\alpha_{m}|^{2} \cdot 
\frac{1}{c} \sum_{x \tpmod{c}} |K(m, x; c)|^{2} \right)^{\frac{1}{2}}.
\end{equation*}
Using the notation~\eqref{eq:K-norm}, we may estimate the diagonal contribution as
\begin{equation*}
\sqrt{cMN} \norm{\alpha}_{2} \norm{\beta}_{2} \norm{K}_{2, \infty}.
\end{equation*}
When $m_{1} \ne m_{2}$, the condition $(m_{1}, m_{2}) = 1$ holds automatically. Hence we may estimate~the off-diagonal contribution as
\begin{equation*}
c^{\frac{1}{4}} \sqrt{N} \norm{\beta}_{2} \mathcal{S}_{\Lambda_{\ell} K}(M, M, c)^{\frac{1}{2}}
\end{equation*}
for some $|\ell| \ll \frac{c}{N}$. To summarise, there holds
\begin{equation*}
\mathcal{S}_{K}(M, N; c) \ll \sqrt{cMN} \norm{\alpha}_{2} \norm{\beta}_{2} \norm{K}_{2, \infty}
 + c^{\frac{1}{4}} \sqrt{N} \norm{\beta}_{2} \mathcal{S}_{\Lambda_{\ell} K}(M, M, c)^{\frac{1}{2}}.
\end{equation*}
We now start iterating the above process deducing
\begin{multline*}
\mathcal{S}_{K}(M, N; c) \ll \sqrt{cMN} \norm{\alpha}_{2} \norm{\beta}_{2} \norm{K}_{2, \infty}
 + \sqrt{cMN} \norm{\alpha}_{2} \norm{\beta}_{2} \norm{\Lambda_{\ell_{1}} K}_{2, \infty}^{\frac{1}{2}}\\
 + c^{\frac{3}{8}} \sqrt{N} M^{\frac{1}{4}} \norm{\alpha}_{2}^{\frac{1}{2}} \norm{\beta}_{2} 
\mathcal{S}_{\Lambda_{\ell_{2}} \Lambda_{\ell_{1}} K}(M, M, c)^{\frac{1}{4}}.
\end{multline*}
Iterating $k-1$ times for sufficiently large $k \in \mathbb{N}$ implies
\begin{multline}\label{iteration}
\mathcal{S}_{K}(M, N; c) \ll \sqrt{cMN} \norm{\alpha}_{2} \norm{\beta}_{2} 
\left(\norm{K}_{2, \infty}+\norm{\Lambda_{\ell_{1}} K}_{2, \infty}^{\frac{1}{2}}
 + \cdots +\norm{\Lambda_{\ell_{k-1}} \cdots \Lambda_{\ell_{1}} K}_{2, \infty}^{\frac{1}{2^{k-1}}} \right)\\
 + c^{\frac{1}{2}-\frac{1}{2^{k+1}}} \sqrt{N} M^{\frac{1}{2}-\frac{1}{2^{k}}} \norm{\alpha}_{2}^{1-\frac{1}{2^{k-1}}} 
\norm{\beta}_{2} \mathcal{S}_{\Lambda_{\ell_{k}} \cdots \Lambda_{\ell_{1}} K}(M, M, c)^{\frac{1}{2^{k}}}.
\end{multline}
To conclude the iteration, we employ the trivial bound
\begin{equation*}
|\mathcal{S}_{K}(M, N; c)| \leq \sqrt{MN} \norm{\alpha}_{2} \norm{\beta}_{2} 
\left(\mathop{\sum \sum}_{\substack{m \sim M, \, n \sim N \\ (m, n) = (mn, c) = 1}} |K(m, n; c)|^{2} \right)^{\frac{1}{2}}.
\end{equation*}
Applying Poisson summation (Lemma~\ref{Fouvry-Kowalski-Michel}) in either of the variables shows
\begin{equation*}
\sum_{n \sim N} |K(m, n; c)|^{2} \ll \frac{N}{c} \sum_{|\ell| \ll \frac{c}{N}} 
\sum_{x \tpmod{c}} |K(m, x; c)|^{2} e \left(-\frac{\ell x}{c} \right) \ll c \norm{K}_{2, \infty}^{2}.
\end{equation*}
Hence
\begin{equation*}
\mathcal{S}_{K}(M, N; c) \ll \min(\sqrt{M}, \sqrt{N}) \sqrt{cMN} \norm{\alpha}_{2} \norm{\beta}_{2} \norm{K}_{2, \infty}.
\end{equation*}
Inserting this estimate into~\eqref{iteration} completes the proof of Theorem~\ref{thm:main}.


\providecommand{\bysame}{\leavevmode\hbox to3em{\hrulefill}\thinspace}
\providecommand{\MR}{\relax\ifhmode\unskip\space\fi MR }
\providecommand{\MRhref}[2]{%
  \href{http://www.ams.org/mathscinet-getitem?mr=#1}{#2}
}
\providecommand{\href}[2]{#2}

\end{document}